\documentclass[11pt]{amsart}

\usepackage[utf8]{inputenc}
\usepackage[T1]{fontenc}
\usepackage[english]{babel}

\usepackage{amsmath, amssymb, amsfonts, amsthm, amsbsy}
\usepackage{dsfont}

\usepackage{graphicx}
\usepackage{tikz}
\usepackage{caption}
\usepackage{subcaption}
\usepackage{float}

\usepackage{enumitem}
\usepackage{relsize}
\usepackage{stackengine}

\usepackage{hyperref}
\usepackage[capitalise,compress]{cleveref}
\usepackage{cite}

\usepackage{color}

\usepackage{comment}
\usepackage{psfrag}

\numberwithin{equation}{section}

\newtheorem{theorem}{Theorem}[section]
\newtheorem{lemma}[theorem]{Lemma}
\newtheorem{proposition}[theorem]{Proposition}
\newtheorem{corollary}[theorem]{Corollary}

\theoremstyle{definition}
\newtheorem{definition}[theorem]{Definition}

\newtheorem{example}[theorem]{Example}

\theoremstyle{remark}
\newtheorem{remark}[theorem]{Remark}

\newcommand{\RR}{\mathbb{R}}

\newcommand{\MM}{\mathbb{M}}

\newcommand{\HH}{\mathbb{H}}
\newcommand{\SPH}{\mathbb{S}}

\newcommand{\Mk}{\mathbb{M}_{\kappa}^d}

\newcommand{\proj}{\operatorname{proj}}

\newcommand{\conv}{\operatorname{conv}}
\newcommand{\pow}{\operatorname{pow}}

\newcommand{\Vold}{\operatorname{Vol}_{\mathbb{M}^d_\kappa}}
\newcommand{\Volk}{\operatorname{Vol}_{\mathbb{M}^k_\kappa}}

\newcommand{\dd}{\mathrm{d}}
\newcommand{\dk}{d_{\kappa}}

\newcommand{\loren}[2]{\langle #1, #2 \rangle_{L}}
\newcommand{\inner}[2]{\langle #1, #2 \rangle_{E}}
\newcommand{\ink}[2]{\langle #1, #2 \rangle_{\kappa}}

\DeclareMathOperator{\arcosh}{arcosh}

\date{\today}

\begin{document}

\sloppy
\title[Laguerre tessellations in spaces of constant curvature]{Laguerre tessellations in spaces of constant curvature and their sectional distributions}

 \author{Panagiotis Spanos}
\address{Panagiotis Spanos, Dipartimento di Matematica,
Sapienza Universit\`a di Roma, Roma, Italy}
\email{panagiotis.spanos@uniroma1.it}

\date{\today}

\renewcommand\footnotemark{}

\keywords{Poisson--Laguerre tessellations, Poisson--Voronoi tessellations, spaces of constant curvature, stochastic geometry, spherical geometry, hyperbolic geometry, sectional distributions}
\thanks{}

\begin{abstract}
We determine the sectional distributions of Poisson--Laguerre tessellations in spaces of non-zero constant curvature, proving closure under totally geodesic sections and deriving explicit weight intensity measures. This includes sections of Poisson--Voronoi tessellations. We construct the \emph{generalised polar region}, a new geometric object which extends existing polar constructions for both positive and negative curvature and provides an alternative representation of Laguerre tessellations. We also give deterministic conditions for tessellation existence and establish sharp integrability criteria for Poisson--Laguerre tessellations.
\end{abstract}

\maketitle

\setcounter{tocdepth}{1}

\baselineskip=16pt

\section{Introduction} 
Poisson--Laguerre tessellations are a natural generalisation of the classical Poisson--Voronoi tessellation. Instead of considering only a point process in the underlying space, each point is equipped with an additional weight, and the corresponding cells are determined by comparing the associated power functions. From a growth-process perspective, the weights may be interpreted as allowing the individual cells to start growing at different times.

Poisson--Laguerre tessellations arise naturally as sections of Poisson--Voronoi tessellations. Let a stationary Poisson--Voronoi tessellation in $\RR^d$ be intersected with a lower-dimensional affine subspace. The resulting collection of cells is again a random tessellation of the subspace, but, as shown by Chiu, van de Weygaert and Stoyan \cite{Weygaert}, it is not itself a Voronoi tessellation. The distribution of such sectional tessellations was identified by Gusakova, Kabluchko and Th\"ale \cite{GKT1}. In particular, they showed that it is a specific Poisson--Laguerre tessellation, which they studied in related work \cite{GKT2,GKT3,GKT4}.

An arbitrary collection of weighted points need not generate a tessellation. Deterministic conditions ensuring the required tessellation properties, together with corresponding conditions for Poisson generators, were developed in the Euclidean setting \cite{Laut1,Laut2}. More recently, Gusakova and in Wolde-L\"ubke \cite{AnnaMathias} studied Poisson--Laguerre tessellations with weight-dependent intensity densities and established that this class is closed under affine sections.

Laguerre tessellations beyond Euclidean space had already appeared in the work of Sugihara \cite{Sugihara}, who introduced the Laguerre diagram on the two-dimensional sphere. In subsequent work, Chaidee and Sugihara \cite{ChaideeSugihara} further studied spherical Laguerre tessellations and established a geometric representation in terms of the central projection of polyhedra in $\RR^3$. Further work on random spherical tessellations includes the study of splitting tessellations \cite{HugThaele} and the typical cell of a spherical Voronoi tessellation \cite{KabThael}. Spherical Voronoi tessellations and other tools from stochastic geometry are used to analyse the ground coverage of satellite constellations (see, for example, \cite{TalgatEtal,DaiEtal}). 

Gusakova and the author \cite{AnnaPan} introduced Laguerre tessellations in hyperbolic space. Using a projective representation of hyperbolic space, they associated to the weighted points a locally polyhedral object, called the polar intersection set, whose facets recover the Laguerre cells under radial projection. This construction allowed them to obtain hyperbolic analogues of several results known in Euclidean space. In particular, they determined the sectional distributions of Poisson--Laguerre tessellations. Hyperbolic Poisson--Laguerre tessellations can also be generated by nuclei on the ideal boundary \cite{DAchilleEtAl2026,AnnaPan}. Since our primary focus is the spherical setting, we restrict attention to nuclei lying inside the underlying space.

The polar construction in \cite{AnnaPan} is closely related to the spherical polyhedral representation considered in \cite{ChaideeSugihara}. However, extending this approach to study sectional properties of arbitrary spherical configurations presents a difficulty. The associated polyhedron may be unbounded, and the radial projections of its facets need not cover the whole sphere. Boundedness of this polyhedron is equivalent to the origin lying in the interior of the convex hull of the generators. One could therefore restrict the construction to configurations satisfying this condition. For Poisson point processes, however, conditioning on this event excludes configurations that occur with positive probability and must be included when determining the sectional distributions of Poisson--Voronoi tessellations.

In this paper, we overcome this obstruction by introducing a new geometric object, the \emph{generalised polar region}, which provides a global representation of spherical Laguerre tessellations without such conditioning. Using this construction, we obtain a unified description of Laguerre tessellations in spaces of non-zero constant curvature and determine the distributions of totally geodesic sections of Poisson--Laguerre tessellations. 
\subsection{Main results}
Let $\Mk$ be a space of constant curvature $\kappa\neq 0$ (see Section~\ref{sec:spaces}). For $h\in (0,\infty)$ and $v,u\in\MM^d_\kappa$ we define the  power function as
$$\pow_\kappa((v,h),u) =h \cos  (\sqrt{\kappa} \dk (v,u) ),
$$
where $\cos ( i x )=\cosh x$, $x\in \RR$. 

Let $A \subset \MM_\kappa^d \times (0,\infty)$
be a countable set of weighted points. 
For $\kappa < 0$, the Laguerre cell
corresponding to 
$(v,h) \in A$ is defined as 
\begin{align}\label{deterministicLaguerre}
C_\kappa((v,h),A) := \{u \in \MM_\kappa^d : \pow_{\kappa}((v,h),u) \leq \pow_{\kappa}((v',h'),u) \text{ for all } (v', h') \in A \}.
\end{align}
For $\kappa > 0$, $A \subset \MM_\kappa^d \times (0,\infty)$,
its corresponding Laguerre cell is
\begin{align}\label{deterministicLaguerre2}
C_\kappa((v,h),A) := \{u \in \MM_\kappa^d : \pow_{\kappa}((v,h),u) \geq \pow_{\kappa}((v',h'),u) \text{ for all } (v', h') \in A \}.
\end{align}
For any curvature $\kappa$, the spatial coordinate $v$ is referred to as the \textit{nucleus} of the cell $C_\kappa((v, h), A)$. 
The collection of all Laguerre cells with non-empty interiors generated by $A$ is referred to as the Laguerre diagram, denoted by
\begin{align}\label{eq:Laguerre_diagrams}
\mathcal{L}^d_{\kappa}(A) := \{C_\kappa((v,h),A) : (v,h)\in A \text{ and } \operatorname{int}\left( C_\kappa((v,h),A)\right) \neq \emptyset \},
\end{align}
where the superscript denotes dimension and we will omit it when there is no ambiguity. We set $\mathcal{L}_\kappa^d(\emptyset)=\emptyset$.

Let $\mu$ be a non-zero, locally finite Borel measure on $(0,\infty)$. We denote by $\mathcal{L}_\kappa^d(\mu)$ the random Laguerre diagram defined for the countable set of points $\eta_\mu$, where $\eta_\mu$ is a Poisson point process with intensity measure $  \operatorname{Vol}_{\Mk} \otimes \, \mu $. 
We say that $\mu$ satisfies the \emph{integrability condition}
$(\mathrm{IC}^d_\kappa)$ if
\begin{equation}\tag{$\mathrm{IC}^d_\kappa$}\label{eq:IC}
\begin{cases}
\mu([1,\infty))<\infty,
& \kappa>0,\\ 
\displaystyle\int_{(0,1]} t^{1-d}\,\mu(\dd t)<\infty,
& \kappa<0,\ d\geq2,\\ 
\displaystyle\int_{(0,1]} \log(1/t)\,\mu(\dd t)<\infty,
& \kappa<0,\ d=1.
\end{cases}
\end{equation}
The integrability condition \eqref{eq:IC} is a natural assumption. As the following theorem shows, it is necessary and sufficient
for $\mathcal{L}_\kappa^d(\mu)$ to be almost surely a tessellation,
conditional on $\eta_\mu\neq\emptyset$. We write $\omega_m$ for the surface area of the unit sphere $\SPH^{m-1}$. Our main results are the following.
\begin{theorem}\label{thrm:iccondition}
Let $d\geq 1$ and $\Mk$ be a space of constant curvature $\kappa\neq 0$. Let $\mu$ be a non-zero, locally finite Borel measure on $(0,\infty)$.
Conditional on $\eta_\mu$ being non-empty, $\mathcal{L}_\kappa(\mu)$ is a tessellation almost surely if and only if $\mu$ satisfies the integrability condition \eqref{eq:IC}. In this case, the tessellation is almost surely face-to-face
and normal, with bounded cells.
\end{theorem}
\begin{theorem}\label{thrm:intersectiondistribution}
Let $1\leq k < d$. Let $\Mk$ be a space of constant curvature $\kappa\neq 0$ and let $L_k$ be a $k$-dimensional totally geodesic subspace of $\Mk$.
Let $\mu$ be a non-zero, locally finite Borel measure on $(0,\infty)$ that satisfies the integrability condition \eqref{eq:IC}. 
Then
$$
\mathcal{L}^d_{\kappa}(\mu)\cap L_k \stackrel{d}{=}\mathcal{L}_{\kappa}^{k}(\tilde{\mu}),
$$
where 
$$
\tilde{\mu}(\dd u)=  \left[ \frac{\omega_{d-k}}{ \vert \kappa\vert^{\frac{d-k}{2}}}  u^k \int_{J_\kappa(u)}      \left\vert t^2-u^2   \right\vert^{\frac{d-k}{2}-1}  t^{1-d}  
\, \mu (\dd t) \right] \dd u,
$$ 
and the integration bounds are $J_\kappa(u) = (0, u)$ if $\kappa < 0$, and $J_\kappa(u) = (u,\infty)$ if $\kappa > 0$.
\end{theorem}
Let $\lambda>0$.
Taking $\mu=\lambda\delta_1$ yields the Poisson--Voronoi
tessellation of intensity $\lambda$.
Theorem~\ref{thrm:intersectiondistribution} therefore gives
the following corollary.
\begin{corollary}\label{thrm:intersectiondistribution2}
Let $1\leq k < d$ and $\lambda>0$. Let $\Mk$ be a space of constant curvature $\kappa\neq 0$ and let $L_k$ be a $k$-dimensional totally geodesic subspace of $\Mk$. Let $\eta_\lambda$ be a Poisson point process on $\Mk$ of intensity $\lambda$. Then
$$
\mathcal{V}^d_{\kappa}(\eta_\lambda)\cap L_k \stackrel{d}{=}\mathcal{L}_{\kappa}^k ({\mu}),
$$ 
where  
$$ 
\mu(\dd t )= \lambda \frac{\omega_{d-k}}{|\kappa|^{\frac{d-k}{2}}} t^k \vert t^2 -  1\vert^{\frac{d-k}{2}-1} \mathbf{1}\{t\in I_\kappa\} \, \dd t ,
$$ 
with $I_\kappa = (0, 1)$ if $\kappa > 0$, and $I_\kappa = (1, \infty)$ if $\kappa < 0$.
\end{corollary}

We also give deterministic conditions under which a collection of weighted points generates a well-defined Laguerre tessellation, as well as conditions under which the usual polar intersection set construction is already sufficient and the generalised polar region is not required.
\section{Preliminaries}
\subsection{Tessellations}
We begin with the deterministic definition of a tessellation.
\begin{definition}
Let $\mathcal{T}=\{ C_i\subset \Mk : i \geq 1 \} $ be an at most countable collection of sets. We call $\mathcal{T}$ a tessellation of $\Mk$ if 
\begin{enumerate}
\item[(i)] Each cell $C\in \mathcal{T}$ is closed and geodesically convex.
\item[(ii)] The union of all cells covers the space $\Mk$.
\item[(iii)] Cells have non-empty interiors and each two different cells have disjoint interior.
\item[(iv)] The tessellation is locally finite.
\end{enumerate}
\end{definition}
We shall mainly consider tessellations with bounded cells.
Every bounded cell is a polytope, that is, a bounded
intersection of finitely many closed geodesic halfspaces.
A non-empty proper face of a polytope is its intersection
with a supporting totally geodesic hyperplane.
A $k$-face is a face of dimension $k$.
As usual, the empty set and the polytope itself are also
regarded as faces.

A tessellation $\mathcal{T}$ is called face-to-face if for every finite collection of cells, their intersection is a face of each individual cell. 
A tessellation $\mathcal{T}$ is called normal if it is face-to-face and, for each $k\in\{0,\ldots,d-1\}$, each $k$-face is contained in the boundary of exactly $d-k+1$ cells.
More detailed exposition can be found \cite[Chapter~10]{SchneiderWeil} for Euclidean tessellations, in \cite[Section~2]{Herold} and \cite[Section~2]{AnnaPan} for hyperbolic tessellations and in \cite[Section~2]{HugThaele} for spherical tessellations. 

\subsection{Spaces of constant curvature}\label{sec:spaces}
We denote by $\MM_\kappa^d$ the standard $d$-dimensional space of constant sectional curvature $\kappa \neq 0$. Depending on the sign of $\kappa$, $\MM_\kappa^d$ corresponds to a scaled sphere ($\kappa > 0$) or a scaled hyperbolic space ($\kappa < 0$).

To treat these spaces in a unified manner, it is convenient to consider their standard embeddings into $\RR^{d+1}$. Let $x,y \in \RR^{d+1}$. We write $\inner{\cdot}{\cdot}$ for the standard Euclidean
inner product and $\Vert\cdot\Vert_2^2$ for its induced squared norm.
Similarly, the Lorentzian inner product and its associated quadratic form are defined as
$$\loren{x}{y}=x_1y_1+\ldots+x_dy_d-x_{d+1}y_{d+1} \quad \text{and} \quad \Vert x \Vert_L^2 = \loren{x}{x}.$$
For $\kappa\neq0$, we denote the curvature-dependent symmetric
bilinear form on $\RR^{d+1}$ by
$$
\ink{x}{y} = 
\begin{cases} 
\inner{x}{y} & \text{if } \kappa > 0, \\
\loren{x}{y} & \text{if } \kappa < 0.
\end{cases}$$

For negative curvature $\kappa<0$, the space $\MM_\kappa^d$ is a scaled hyperbolic space. We model $\MM_\kappa^d$ using the upper sheet of the scaled hyperboloid
$\MM_\kappa^d = \{ x\in \RR^{d+1}: \Vert x \Vert_L^2 = \kappa^{-1} , \, x_{d+1}>0 \} \subset \RR^{d+1}$.
The intrinsic geodesic distance $\dk(x,y)$ between $x, y \in \MM_\kappa^d$ is given by
$$\dk (x,y)= \frac{1}{\sqrt{-\kappa}} \arcosh 
( \kappa \loren{x}{y} ).
$$

For positive curvature $\kappa > 0$, the space $\MM_\kappa^d$ is modeled as a sphere of radius $1/\sqrt{\kappa}$, defined by
$\MM_\kappa^d = \left\{ x\in\RR^{d+1}: \Vert x \Vert_2^2 =\kappa^{-1}\right\} \subset \RR^{d+1}$.
The intrinsic geodesic distance in this case is given by the spherical metric,
$$ \dk (x,y) =  \frac{1}{ \sqrt{ \kappa } } \arccos ( \kappa \inner{x}{y} ).
$$

We define the \emph{generalised cosine function} 
$$\cos_\kappa (r )=\begin{cases}
\cos(\sqrt{\kappa}\,r), & \text{if } \kappa > 0, \\ 
\cosh(\sqrt{-\kappa}\,r), & \text{if } \kappa < 0
\end{cases}
$$
and the \emph{generalised sine function} 
$$\sin_\kappa( r )=
\begin{cases}
\dfrac{1}{\sqrt{\kappa}}\sin(\sqrt{\kappa}\,r), & \text{if } \kappa > 0, \\
\dfrac{1}{\sqrt{-\kappa}}\sinh(\sqrt{-\kappa}\,r), & \text{if } \kappa < 0.
\end{cases}
$$
These functions satisfy the identity
\begin{equation}\label{eq:genpythagorean}
\cos_\kappa^2(r)+\kappa \sin_\kappa^2(r)=1.
\end{equation}
The  power functions can be expressed in a unified way for $\kappa\neq 0$ as
$$\pow_\kappa((v,h),u) =h \cos_\kappa ( \dk (v,u) )=\kappa\ink{hv}{u},
$$
for $h\in (0,\infty)$ and $v,u\in\MM^d_\kappa$.

Another crucial notion for our construction is the radial projection from the ambient space onto the manifold $\MM_\kappa^d$. We will denote it by $\proj_{\MM_\kappa^d}$. For $\kappa>0$, the radial projection is 
$$\proj_{\MM_\kappa^d}(x) = \frac{1}{ \sqrt{\kappa \Vert x \Vert_2^2} } x$$
and it is well-defined for all non-zero vectors in the ambient space, $\Omega_\kappa = \RR^{d+1} \setminus \{0\}$.
For negative curvature ($\kappa < 0$), the projection is 
$$\proj_{\MM_\kappa^d}(x) = \frac{1}{ \sqrt{ \kappa \Vert x \Vert_L^2 } } x$$ 
and it is well-defined for elements within the interior of the upper light cone, $\Omega_\kappa = \{x \in \RR^{d+1} : \Vert x \Vert_L^2 < 0, \, x_{d+1} > 0\}$.

\section{Laguerre Tessellations}
In this section, we study the geometric properties of Laguerre
cells and give deterministic conditions under which they form
a tessellation.
\subsection{Laguerre tessellations}
Let $A\subset\MM_\kappa^d\times(0,\infty)$ be countable.
For each $(v,h)\in A$, the ambient representation of the
power function gives
\[
C_\kappa((v,h),A)
=
\left\{
u\in\MM_\kappa^d:
\ink{hv-h'v'}{u}\geq0
\text{ for all }(v',h')\in A
\right\}.
\]
Thus each Laguerre cell is the intersection of $\MM_\kappa^d$
with a possibly infinite family of closed ambient halfspaces,
and is therefore closed and geodesically convex.

While the Laguerre diagram is well-defined for any arbitrary set of weighted points, it only forms a well-defined tessellation, with bounded cells, under certain regularity and general position conditions. Let $A \subset \MM^d_\kappa\times (0,\infty)$ be an at most countable set of weighted points. In accordance with the regularity conditions for $\RR^d$ and $\HH^d$ in \cite{Laut1,AnnaPan,GKT3}, we impose the following regularity conditions
\begin{enumerate} 
\item[(P1)] 
Let $\kappa\neq 0$, for every $u\in \MM^d_\kappa$ there exists $s\in\RR$, such that
$$ 0< \#\bigl\{(v,h)\in A:\langle hv,u\rangle_\kappa>s\bigr\} <\infty.$$
\end{enumerate}
\begin{remark}
Condition \textup{(P1)} ensures that the extremal power at
each point is obtained among finitely many generators.
On the sphere, it allows infinitely many generators with
weights tending to zero, as may occur for Poisson--Laguerre
models satisfying \eqref{eq:IC}.
For negative curvature, it is a variant of condition
\textup{(P1)} in \cite{AnnaPan}.
However, on the sphere, requiring finiteness only above
positive thresholds would not suffice. Indeed, the configuration
$A=\{(v,1/n):n\geq1\}$ satisfies this requirement, but its
powers at the antipodal point $-v$ are $-1/n$, whose supremum
is not attained. Thus $-v$ belongs to no cell.
The additional non-emptiness requirement in \textup{(P1)}
excludes this pathology without forcing $A$ itself to be finite.
\end{remark}
To ensure that all cells in the tessellation are bounded, we require a curvature-dependent condition. We define condition \textup{(P2)} according to curvature.
\begin{enumerate}
\item[(P2)] 
\begin{itemize}
	\item For $\kappa<0$: For every $x \in \{ x \in \mathbb{R}^{d+1} : x_{d+1}=1 \text{ and } \loren{x}{x}=0\}$ and every $s >0 $, there exists $(v,h) \in A$ such that
$$ -h\loren{v}{x} < s.$$
	\item For $\kappa>0$: The set $A$ is non-empty.
	\end{itemize}
\end{enumerate} 
\begin{remark}
For comparison, the Euclidean analogue of \textup{(P1)}
requires that, for every $u\in\RR^d$ and $s\in\RR$, only
finitely many $(v,h)\in A$ satisfy
$\Vert v-u\Vert_2^2+\log h<s$.
The Euclidean analogue of \textup{(P2)} is
$\conv\{v:(v,h)\in A\}=\RR^d$.
While \textup{(P1)} imposes a finiteness requirement on the
weighted generators, \textup{(P2)} expresses a global geometric
condition whose formulation depends on the curvature.
\end{remark}
Furthermore, we impose the general position condition.
\begin{enumerate}
\item[(P3)]
For every $m\in\{2,\ldots,d+2\}$, there is no
$(d-m+2)$-dimensional totally geodesic subspace of
$\MM_\kappa^d$ on which the power functions of $m$
distinct elements of $A$ coincide.
\end{enumerate}
With these conditions established, we can state the following proposition, establishing tessellations.
\begin{proposition}\label{prop:technicalprop}
Let $\Mk$ be a space of constant curvature $\kappa\neq 0$. Then
\begin{enumerate}
\item If $A\neq\emptyset$ and $A$ satisfies conditions \textup{(P1)}, then the Laguerre diagram $\mathcal{L}_\kappa(A)$ forms a face-to-face tessellation of $\MM^d_\kappa$.
\item If $A$ satisfies conditions \textup{(P1)} and \textup{(P2)}, then the Laguerre diagram $\mathcal{L}_\kappa(A)$ forms a face-to-face tessellation of $\Mk$ with bounded cells.
\item If, in addition to \textup{(P1)} and \textup{(P2)}, $A$ satisfies \textup{(P3)}, then the Laguerre diagram $\mathcal{L}_\kappa(A)$ is a normal tessellation of $\MM^d_\kappa$.
\end{enumerate}
\end{proposition}
\begin{proof}
The proof follows the same arguments as in
\cite[Proposition~2.2.5 and Corollary~2.2.7]{Laut1}
and \cite[Theorem~3.8]{AnnaPan}.
The only adjustment concerns \textup{(P1)}, which here
requires a finite, non-empty set of competitors at one
suitable threshold for each point, rather than finiteness
at every threshold. 

We write $f_p(v):=\pow_\kappa (p,v)$ for every $v\in\Mk$ and $p\in A$. 
For $\kappa<0$, we fix $u\in\Mk$. By \textup{(P1)}, there exists $b>0$ such that the set $I_u:=\{ p\in A: f_p(u)<b\}$ is finite and non-empty. Therefore we can define $m_u:=\min_{p\in I_u} f_p(u)<b$.
The triangle inequality gives 
$$
e^{-r\sqrt{-\kappa}} f_p(u) \leq f_p(v) \leq e^{r\sqrt{-\kappa}} f_p(u),
$$
for every $v\in\Mk$ that satisfies $d_\kappa(v,u)\leq r$. Because $m_u<b$, we can choose $r>0$ so that $e^{r\sqrt{-\kappa}} m_u < e^{-r\sqrt{-\kappa}} b$, which shows that no generator outside $I_u$ can determine a cell in the $r$-ball around $u$.

For $\kappa>0$, we apply \textup{(P1)} at each standard basis vector 
$e_j/\sqrt\kappa$, and their antipodals $-e_j/\sqrt\kappa$, for $j=1,\ldots, d+1$. This shows that all the vectors $hv$ are bounded,
in particular $\Vert hv \Vert_2\le C$, uniformly over $(v,h)\in A$.
Fix $u\in\Mk$ and choose $s$ such that
$I_u:=\{p\in A:f_p(u)>s\}$
is finite and non-empty. Its maximal value $M_u=\max_{p\in I_u}f_p(u)$ exceeds $s$.
The inequality 
$$|f_p(v)-f_p(u)|=\kappa |\inner{hv}{v-u}|\le \kappa C\|v-u\|_2$$
shows that, for $v$ sufficiently close to $u$ so that $\|v-u\|_2<(M_u-s)/(2\kappa C)$, 
every generator outside $I_u$ has strictly smaller score than
a maximizing generator in $I_u$.
\end{proof}

\section{Generalised polar region}
In this section, we give a geometric representation of Laguerre
tessellations in spaces of non-zero constant curvature, extending
the polar constructions in \cite{AnnaPan,Sugihara}.

We first describe the geometric idea in the unweighted case,
which corresponds to the Voronoi tessellation. At each nucleus,
consider the tangent hyperplane to $\Mk$ and take the closed
halfspace containing $\Mk$. The intersection of these halfspaces
is the \emph{polar intersection set}. The idea is to recover
the Voronoi cells by radially projecting its facets onto $\Mk$.

\begin{figure}[htbp]
    \centering

    \begin{subfigure}{0.32\textwidth}
        \centering
        \includegraphics[width=\linewidth]{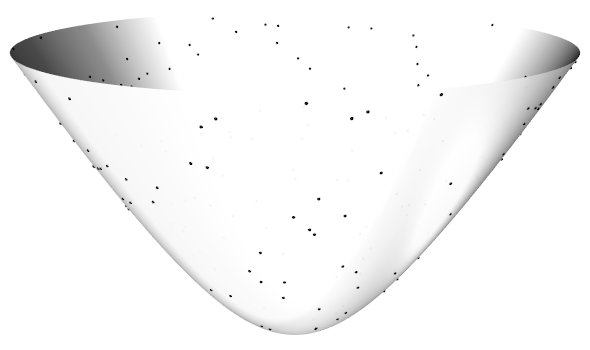}
        \caption{Points in $\HH^2$.}
        \label{fig:image1}
    \end{subfigure}
    \hfill
    \begin{subfigure}{0.32\textwidth}
        \centering
        \includegraphics[width=\linewidth]{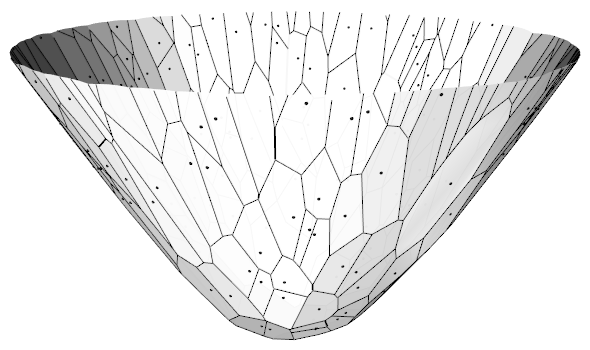}
        \caption{polar intersection set in $\RR^3$.}
        \label{fig:image2}
    \end{subfigure}
    \hfill
    \begin{subfigure}{0.32\textwidth}
        \centering
        \includegraphics[width=\linewidth]{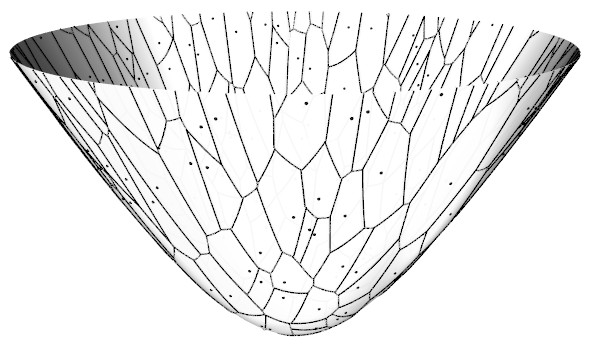}
        \caption{Voronoi cells in $\HH^2$.}
        \label{fig:image3}
    \end{subfigure}

    \caption{The construction of the Voronoi tessellation in $\HH^2$, via radial projection of the polar intersection set of the points.}
    \label{fig:three_images}
\end{figure}

On the sphere, however, this set may be unbounded, and some
positive rays from the origin may never meet its boundary.
The projected facets then fail to cover the whole sphere.
We overcome this obstruction by introducing the
\emph{generalised polar region}.

To generalise this construction to Laguerre tessellations,
where each nucleus carries a weight, we radially scale the
nuclei in the ambient space. The corresponding hyperplanes
are shifted parallel to their original tangent positions.
The choice of scaling factor in the hyperbolic setting
follows from the projective construction developed in
\cite{AnnaPan}. We extend this construction to spaces of
non-zero constant curvature and show that it provides
an alternative geometric representation of Laguerre
tessellations.

Let $A \subset \MM_\kappa^d \times (0,\infty)$  be an at most countable set of weighted points. We map each weighted point into the ambient space $\RR^{d+1}$ via the radial scaling function $\alpha_\kappa : \MM_\kappa^d \times (0,\infty) \to \RR^{d+1}$  defined by  $$\alpha_\kappa(v,h) = h v. $$

By the ambient representation of the power function, we have 
$\pow_\kappa((v,h),u) = \kappa \ink{\alpha_\kappa(v,h)}{u}  $.
In the hyperbolic case ($\kappa < 0$), the Laguerre cell consists of $u \in \MM_\kappa^d $ satisfying
$$ \pow_{\kappa}((v,h),u) \leq \pow_{\kappa}((v',h'),u),\quad \text{ for all } (v', h') \in A.
$$
Dividing by the negative constant $\kappa$ reverses the inequality, yielding 
$$
  \ink{\alpha_\kappa(v,h)}{u}   \geq \ink{\alpha_\kappa(v',h')}{u},\quad  \text{ for all } (v', h') \in A .
$$ 
In the spherical case ($\kappa > 0$), the cell is defined by maximizing the power function. Dividing by the positive constant $\kappa$ preserves the inequality, yielding the same relation.
Thus, for any curvature $\kappa \neq 0$, a point $u \in C_\kappa((v,h),A)$ if and only if
$$\ink{\alpha_\kappa(v,h)}{u}   \geq \ink{\alpha_\kappa(v',h')}{u}, \quad \text{ for all } (v',h') \in A. $$

For a mapped point $x = \alpha_\kappa(v,h)$, we define the \emph{polar hyperplane} $x^{\circ}_\kappa = \left\{ y \in \RR^{d+1} : \ink{x}{y} = \kappa^{-1} \right\}$. This represents the translation of the tangent hyperplane at $v$ according to the weight $h$. The corresponding closed half-space is given by
$$
(x^{\circ}_{\kappa})^+= \left\{ y \in \RR^{d+1} : \ink{x}{y} \leq \kappa^{-1} \right\}.
$$
We will omit the subscript $\kappa$ when the curvature is obvious. We define the \emph{polar intersection set} $P(A) \subset \RR^{d+1}$ as the intersection of all such half-spaces
$$
P(A) = \bigcap_{(v,h) \in A} (\alpha_\kappa(v,h)^{\circ})^+ .
$$	
For each $(v,h)\in A$, define the corresponding
possibly empty face by
$$
F((v,h),A):=\alpha_\kappa(v,h)^\circ\cap P(A).
$$

The polar construction has a projective interpretation
(see \cite{AnnaPan}). In negative curvature, the hyperboloid
model admits a global projective representation.
On the sphere, the analogous projectivization identifies
antipodal points and is therefore not injective.
Its restriction to an open hemisphere, however, identifies
that hemisphere with an affine chart of projective space,
so the construction can be carried out locally.

The difficulty is that the polar intersection set constructed
above may be unbounded, so there may be directions in which
the corresponding positive ray does not intersect its boundary.
We first consider configurations for which this obstruction
does not occur. More precisely, we impose the condition
\begin{itemize}
\item[(P2$'$)]
For $\kappa>0$,
$\max_{(v,h)\in A}\inner{hv}{u}>0$
for every $u\in\Mk$.
\end{itemize}
The weights do not affect the sign, since they are positive.
Therefore, condition \textup{(P2$'$)} is implies
that, for every $u\in\Mk$, there exists $(v,h)\in A$ such that
$\inner{v}{u}>0$. Geometrically, this means that no open
hemisphere is empty of nuclei.

\begin{figure}[htbp]
    \centering

    \begin{subfigure}{0.32\textwidth}
        \centering
        \includegraphics[width=\linewidth]{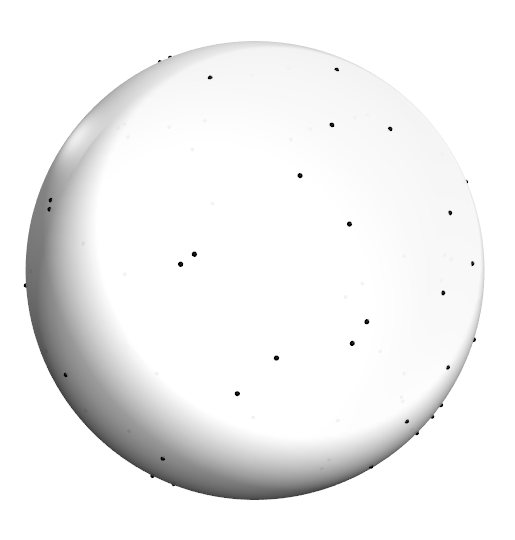}
        \caption{Points in $\mathbb{S}^2$.}
        \label{fig:image4}
    \end{subfigure}
    \hfill
    \begin{subfigure}{0.32\textwidth}
        \centering
        \includegraphics[width=\linewidth]{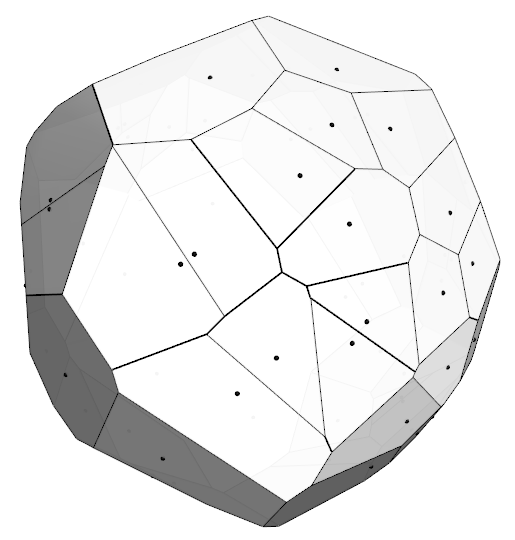}
        \caption{polar intersection set in $\RR^3$.}
        \label{fig:image5}
    \end{subfigure}
    \hfill
    \begin{subfigure}{0.32\textwidth}
        \centering
        \includegraphics[width=\linewidth]{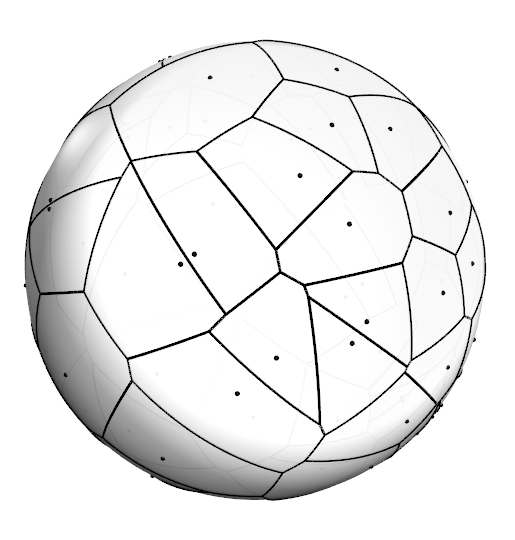}
        \caption{Voronoi cells in $\mathbb{S}^2$.}
        \label{fig:image6}
    \end{subfigure}

    \caption{The construction of the Voronoi tessellation in $\mathbb{S}^2$, via radial projection of the polar intersection set of the points.}
    \label{fig:three_images2}
\end{figure}

The following proposition establishes that the radial projection of these faces onto the manifold exactly recovers the Laguerre cells.  
\begin{proposition}\label{prop:project1}
Let $\kappa\neq0$ and let
$A\subset\MM_\kappa^d\times (0,\infty)$ be an at most countable set of weighted points.
If $\kappa>0$, assume additionally that $A$ satisfies (P2$'$).
Then for every $(v,h)\in A$,
$$
C_\kappa((v,h), A) = \proj_{\MM_\kappa^d} \big( F((v,h), A) \cap \Omega_\kappa \big).
$$
\end{proposition}
\begin{proof}
Let $u \in C_\kappa((v,h), A)$. Then we have $\ink{\alpha_\kappa(v,h)}{u} \geq \ink{\alpha_\kappa(v',h')}{u}$ for all $(v',h') \in A$. We consider the ray generated by $u$ in $\RR^{d+1}$ by defining $y_t = t u$ for $t > 0$. There is a $y_{t_0}$ such that it lies on the bounding hyperplane $\alpha_\kappa(v,h)^{\circ}$. This requires $\ink{\alpha_\kappa(v,h)}{t_0 u} = \kappa^{-1}$, hence
$$
 t_0 = \frac{\kappa^{-1}}{\ink{\alpha_\kappa(v,h)}{u}}= \frac{1}{\pow_\kappa((v,h),u)},
$$
and $t_0>0$. If $\kappa<0$, this follows from the
strict positivity of the hyperbolic power function. If $\kappa>0$, condition
\textup{(P2$'$)} guarantees that the power function obtains its maximum.

Since $t_0 > 0$, multiplying the cell inequality by $t_0$ preserves the inequality. Therefore
$
t_0 \ink{\alpha_\kappa(v,h)}{u} \geq t_0 \ink{\alpha_\kappa(v',h')}{u}$ and hence  $\kappa^{-1} \geq \ink{\alpha_\kappa(v',h')}{y_{t_0}}$ for all $(v',h') \in A$.

This implies that $y_{t_0} \in (\alpha_\kappa(v',h')^{\circ})^{+}$ for every point in $A$, and thus $y_{t_0} \in P(A)$. Because $y_{t_0}$ also satisfies equality for the hyperplane of $(v,h)$, it belongs to the face $F((v,h), A)$. Thus $y_{t_0}\in F((v,h),A)\cap\Omega_\kappa$,
and its radial projection is $u$.

Conversely, let $y\in F((v,h),A)\cap\Omega_\kappa$. Then by definition
\begin{equation}\label{eq:inequalityF}
\ink{\alpha_\kappa(v,h)}{y}\geq \ink{\alpha_\kappa(v',h')}{y},
\end{equation}
for all $(v',h')\in A$. Denote $u=\proj_{\Mk}(y)$ the projection to $\Mk$ and $t_u$ the unique positive number such that $y=t_u u$. Multiplying both sides of inequality (\ref{eq:inequalityF}) by $t_u^{-1}$ we obtain
$$
\ink{\alpha_\kappa(v,h)}{u}\geq \ink{\alpha_\kappa(v',h')}{u},
$$
for all $(v',h')\in A$, which implies the desired equivalence.
\end{proof}

\subsection{Generalised polar region}
In this subsection we restrict attention to $\kappa>0$,
since the hyperbolic case is already covered by
Proposition~\ref{prop:project1}. 

To treat configurations that fail \textup{(P2$'$)}, we
introduce an additional polar intersection set.
Failure of \textup{(P2$'$)} means that some open hemisphere
contains no nuclei. The geometric idea is to supplement
the original projective charts with charts centred at the
antipodal points of the nuclei. This leads to a second
polar intersection set. A further step, involving closures,
will be needed to recover the whole tessellation.

For every $u\in \Mk$, we define the function 
$$
m(u)=\sup_{(v,h)\in A} \inner{hv}{u}.
$$
To account for directions in which $m(u)<0$, we associate
to each mapped point $x=\alpha_\kappa(v,h)$ the closed half-space given by
$$
(x_\kappa^{\circ})^-=\left\lbrace y\in\RR^{d+1}: \inner{x}{y}\leq - \kappa^{-1}\right\rbrace.
$$
Correspondingly we define the negative polar intersection set for the negative values of the inner product
$$
P^{-}(A)=\bigcap_{(v,h)\in A} (\alpha_\kappa(v,h)^{\circ})^-.
$$
Note that this is a closed and convex set and, whenever it is non-empty, it is always unbounded. The set $P^{-}(A)$ never contains the origin of $\RR^{d+1}$. We write $P^+(A):=P(A)$, to distinguish the two sets.

We define the \emph{generalised polar region} as
$$
PR(A):= P^+(A) \setminus \mathrm{int} (P^- (A) ).
$$
If $P^-(A)=\emptyset$, then $PR(A)=P(A)$.
In particular, this holds when \textup{(P2$'$)} is satisfied.
The set $P^{-}(A)$ is non-empty if and only if there exists $u\in\Mk$ such that $m(u)<0$. Thus, if $P^-(A)$ is empty, then no second family of projective charts is needed. the original charts already represent
the tessellation, with directions satisfying $m(u)=0$,
if any, appearing as boundary directions of these charts.

\begin{example}
Consider the unweighted case corresponding to the ordinary Voronoi tessellation. For every nucleus $v\in\Mk$, the positive polar hyperplane is the tangent hyperplane to $\Mk$ at $v$. The corresponding negative hyperplane is the parallel tangent hyperplane at the antipodal point $-v$. We then consider the half-space lying beyond this hyperplane and not containing the origin. The intersection of these half-spaces over all nuclei defines the negative polar intersection set $P^-(A)$.
\end{example}

We will show that the cells of the tessellation can be
recovered from the two polar intersection sets $P^+(A)$
and $P^-(A)$. For each $(v,h)\in A$, define the possibly
empty faces
$$
F^+((v,h),A):= P^+(A) \cap  \alpha_\kappa(v,h)^{\circ} \quad
\text{and}\quad F^-((v,h),A):=P^-(A) \cap (-\alpha_\kappa(v,h))^{\circ}.
$$
Note that as defined in the previous section $P(A)=P^+(A)$ .
The non-empty faces form boundary pieces of $PR(A)$.
The following theorem describes how these faces recover
the Laguerre cells through radial projection.
\begin{theorem}\label{prop:cellsprojection}
Let $\kappa>0$ and let $A\subset \Mk\times (0,\infty)$ be an at most countable set of weighted points. Fix $(v,h)\in A$ and assume that $C_\kappa((v,h),A)$ is full-dimensional. Then the Laguerre cell satisfies
$$
C_\kappa((v,h),A)= \mathrm{cl}\left( 
\proj_{\Mk}(F^+((v,h),A)) \cup \proj_{\Mk}(F^-((v,h),A))
\right),
$$
where $\mathrm{cl}$ denotes closure in $\Mk$.
\end{theorem}
\begin{proof}
Let $u\in C_\kappa((v,h),A)$ and suppose that $m(u)>0$. 
The argument of Proposition~\ref{prop:project1} shows that
$u\in\proj_{\Mk}(F^+((v,h),A))$.
Conversely, the same argument shows that if
$u\in\proj_{\Mk}(F^+((v,h),A))$, then
$u\in C_\kappa((v,h),A)$ and $m(u)>0$.

Now let $u\in C_\kappa((v,h),A)$ and suppose that $m(u)<0$.
Then $m(u)=\inner{\alpha_\kappa(v,h)}{u} \geq  \inner{\alpha_\kappa(v',h')}{u} $ for every $(v',h')\in A$.
As in the positive case, the ray $y_t=tu$ meets the hyperplane 
$-(\alpha_\kappa(v,h))^\circ$
at the point $y_{t_0}=t_0u$ that satisfies $\inner{\alpha_\kappa(v,h)}{t_0 u }=-\kappa^{-1}$. Hence $t_0=-1/(\kappa m(u))>0$.
Multiplying the cell inequalities by $t_0$ gives
$y_{t_0}\in(\alpha_\kappa(v',h')^\circ)^-$
for every $(v',h')\in A$, and thus $y_{t_0}\in P^-(A)$.
Since $y_{t_0}$ also lies on $(-\alpha_\kappa(v,h))^\circ$,
it belongs to $F^-((v,h),A)$, and its radial projection is $u$.

Conversely, let $y\in F^-((v,h),A)$ and set $u=\proj_{\Mk}(y)$. 
Then $y=tu$ for some $t>0$. By the definition of $F^-((v,h),A)$,
$$ \inner{\alpha_\kappa(v',h')}{y} \leq
\inner{\alpha_\kappa(v,h)}{y} = -\kappa^{-1}
$$
for every $(v',h')\in A$. Since $t>0$, dividing by $t$ gives
$u\in C_\kappa((v,h),A)$ and
$m(u)=-\kappa^{-1}/t<0$.

So far we have proven that 
$$
\proj_{\Mk} ( F^{-} ( (v,h) , A)) =
C_\kappa ( (v,h), A) \cap \lbrace u : m(u) < 0 \rbrace,
$$
and the analogous equality in the positive case. Since $C_\kappa((v,h),A)$ is closed, it follows that
$$
 \mathrm{cl}\left( 
\proj_{\Mk}(F^+((v,h),A)) \cup \proj_{\Mk}(F^-((v,h),A))
\right) \subseteq C_\kappa((v,h),A).
$$

Now let $u\in C_\kappa((v,h),A)$ satisfy $m(u)=0$. 
Since the cell is full-dimensional, it cannot be contained
in the great subsphere orthogonal to $\alpha_\kappa(v,h)$.
Thus there exists $w\in C_\kappa((v,h),A)$ such that
$m(w)=\inner{\alpha_\kappa(v,h)}{w}\neq 0$.

Let $\epsilon \in (0,1)$, and consider $x_\epsilon : = (1-\epsilon)u+ \epsilon w$, this is a line in $\RR^{d+1}$, but its projection
$\proj_{\Mk}(x_\epsilon)\in C_\kappa((v,h),A)$. 
Since
$\inner{\alpha_\kappa(v,h)}{x_\epsilon}
=\epsilon m(w)\neq0$,
we have $x_\epsilon\neq0$. We may therefore define
$u_\epsilon:=\proj_{\Mk}(x_\epsilon)
=t_\epsilon x_\epsilon$ for some $t_\epsilon>0$.
The cell inequalities are preserved by convex combinations
and positive rescaling, so
$u_\epsilon\in C_\kappa((v,h),A)$. Consequently,
$m(u_\epsilon)=\epsilon t_\epsilon m(w)\neq0$.

Thus for all $1>\epsilon>0$, either $u_\epsilon \in \proj_{\Mk} ( F^{-} ( (v,h) , A))$ or 
$u_\epsilon \in \proj_{\Mk} ( F^{+} ( (v,h) , A))$. As $\epsilon\downarrow0$, we have $x_\epsilon\to u$,
and continuity of radial projection gives $u_\epsilon\to u$.
We conclude 
$$
u\in \mathrm{cl}\left( 
\proj_{\Mk}(F^+((v,h),A)) \cup \proj_{\Mk}(F^-((v,h),A))
\right).
$$
\end{proof}
The full-dimensionality assumption imposes no additional
restriction on the cells of $\mathcal L_\kappa(A)$,
which by definition have non-empty interiors. 
\begin{remark}
Let $C_\kappa((v,h),A)\in\mathcal L_\kappa(A)$.
If the negative face $F^-((v,h),A)$ is non-empty and bounded, then the corresponding Laguerre cell cannot contain its nucleus $v$. Indeed, if 
$v\in C_\kappa((v,h),A)$, then $\inner{\alpha_\kappa(v,h)}{v}>0$,
whereas the existence of $F^-((v,h),A)$ implies that the cell contains points $u$ for which $\inner{\alpha_\kappa(v,h)}{u}<0$.
By convexity of the cell, it must therefore contain a direction $w$ satisfying $\inner{\alpha_\kappa(v,h)}{w}=0$.
Consequently, for any $y\in F^-((v,h),A)$ and $t\geq0$,
we have $y+tw\in F^-((v,h),A)$.
Since $w\neq0$, this shows that the negative face is unbounded,
a contradiction. 
\end{remark}

\begin{remark}
Suppose that $A$ is finite and contains $d+2$ weighted points whose ambient
vectors are affinely independent. Their convex hull has
non-empty interior in $\RR^{d+1}$. Choose a point $z$ in
this interior, with $z\notin\alpha_\kappa(A)$, and express
the weighted vectors relative to $z$, replacing each $hv$
by $hv-z$.

This translation leaves the Laguerre cells unchanged, since
$\inner{hv-z}{u}\geq\inner{h'v'-z}{u}$ if and only if $\inner{hv}{u}\geq\inner{h'v'}{u}$.
Each translated vector is non-zero and can again be represented
by a nucleus on $\Mk$ and a positive weight. Moreover, the
origin lies in the interior of their convex hull, so the
translated generators satisfy \textup{(P2$'$)}. The usual
polar construction therefore provides a global representation
of the same tessellation.

This gives an alternative construction under the stated
assumption, but does not cover all configurations considered
here. For example, $A$ may contain only one weighted point, or be infinite.
\end{remark}

\section{Poisson--Laguerre tessellations}
In this section we prove Theorem~\ref{thrm:iccondition},
which characterises the measures generating Poisson--Laguerre
tessellations and establishes their regularity.

Let $\mu$ be a non-zero, locally finite Borel measure on
$(0,\infty)$. Recall that $\eta_\mu$ is a Poisson point process
on $\Mk\times(0,\infty)$ with intensity measure
$$
\Lambda_\mu =\Vold \otimes \,\mu.
$$
The intensity is homogeneous in the spatial coordinate,
while $\mu$ determines the intensity in the weight coordinate.
We first consider positive curvature.
\begin{lemma}\label{le:firsttechnical}
Let $\kappa>0$ and let $\mu$ be a non-zero, locally finite
Borel measure on $(0,\infty)$.
Conditional on $\eta_\mu\neq\emptyset$, the process
$\eta_\mu$ satisfies \textup{(P1)} almost surely if and
only if $\mu([1,\infty))<\infty$.
\end{lemma}
\begin{proof}
Suppose first that $\mu([1,\infty))<\infty$.
By local finiteness, $\mu([\varepsilon,\infty))<\infty$
for every $\varepsilon>0$. Since $\Mk$ has finite volume,
almost surely there are only finitely many generators
with weight at least $\varepsilon$, simultaneously for
every $\varepsilon>0$.

If $\mu((0,\infty))<\infty$, then $\eta_\mu$ is finite
almost surely and, conditional on being non-empty,
satisfies \textup{(P1)}. If $\mu((0,\infty))=\infty$, then
every non-empty open subset of $\Mk$ contains
infinitely many nuclei almost surely. Using a countable
basis, it follows that the nuclei are almost surely dense.
Consequently, for every $u\in\Mk$, there exists
$(v,h)\in\eta_\mu$ such that $\inner{hv}{u}>0$.
Choose $s$ with $0<s<\inner{hv}{u}$. The set
$\{(v',h')\in\eta_\mu:\inner{h'v'}{u}>s\}$
is non-empty and finite. Indeed,
$\inner{h'v'}{u}\leq h'/\kappa$, so every such generator
has $h'>\kappa s$. Thus \textup{(P1)} holds.

Conversely, suppose that $\mu([1,\infty))=\infty$.
Local finiteness implies that $\mu((a,\infty))=\infty$
for every $a>0$.
Fix $u\in\Mk$ and $0<c<1/\kappa$, then the non-empty open cap $U_c=\{v\in\Mk:\inner{v}{u}>c\}$
For each integer $n\geq1$,
$\Lambda_\mu\bigl(U_c \times(n/c,\infty)\bigr)=\infty$.
Hence, almost surely, infinitely many points of the point process satisfy
$\inner{hv}{u}>n$, simultaneously for every $n\geq1$.
Therefore there is no threshold above which the number
of generators is finite and non-zero, so
\textup{(P1)} fails almost surely.
\end{proof}
\begin{lemma}
Let $\kappa>0$ and let $\mu$ be a non-zero, locally finite Borel measure on $(0,\infty)$. Then $\eta_\mu$ satisfies \textup{(P3)} almost surely.
\end{lemma}
\begin{proof}
Since $\mu$ is locally finite and $\Mk$ is compact,
the restriction of $\eta_\mu$ to $\Mk\times[1/n,n]$
is almost surely finite for every $n\geq1$.
Conditional on the number of points and their weights,
the vectors $hv$ are independent and uniform on Euclidean
spheres of radii $h/\sqrt\kappa$. Every proper affine
subspace intersects each such sphere in a set of surface
measure zero. Successive conditioning therefore shows that
any at most $d+2$ of these vectors are affinely independent
almost surely. Taking a countable intersection over $n$
establishes this for the entire process.

For $m$ distinct ambient vectors $a_i=h_iv_i$, equality
of their powers at $u$ means $\inner{a_i-a_1}{u}=0$, for every $i=2,\ldots,m$.
Affine independence gives $m-1$ independent equations.
Their common solution set on $\Mk$ has dimension
$d-m+1$ for $2\leq m\leq d+1$ and is empty for $m=d+2$.
Thus \textup{(P3)} holds almost surely.
\end{proof}
\begin{proof}[Proof of Theorem~\ref{thrm:iccondition}]
For $\kappa>0$, suppose that $\mu$ satisfies \eqref{eq:IC}.
Conditional on $\eta_\mu\neq\emptyset$, the preceding
lemmata establish \textup{(P1)} and \textup{(P3)} almost
surely, while \textup{(P2)} holds automatically.
Proposition~\ref{prop:technicalprop} therefore
implies that the diagram is almost surely a face-to-face,
normal tessellation with bounded cells.

Conversely, if $\mu([1,\infty))=\infty$, the argument
in Lemma~\ref{le:firsttechnical} shows that, at any fixed
$u\in\Mk$, the power functions are almost surely
unbounded above. Thus $u$ belongs to no Laguerre cell,
and the diagram is almost surely not a tessellation.

For $\kappa<0$, the result follows from
\cite[Theorem~1.1 and Proposition~4.5]{AnnaPan}
by rescaling the hyperbolic model to curvature $-1$.
This rescaling multiplies the spatial intensity by a
positive constant, leaving the integrability condition
unchanged, and preserves the tessellation properties
and boundedness of the cells.
\end{proof}

\section{Proof of Theorem~\ref{thrm:intersectiondistribution}}
In this section we will use the polar intersection set and the generalised polar region to study the distribution of the intersection of a
Poisson--Laguerre tessellation with a $k$-dimensional totally geodesic subspace of $\Mk$.
In general, the Poisson--Laguerre tessellations are closed under intersections, i.e. the intersection of a Poisson--Laguerre tessellation in $\Mk$ with a totally geodesic subspace of $\Mk$ is itself a Poisson--Laguerre tessellation.

Before we move to the proof of the theorem, we describe how the subspaces can be expressed in $\Mk$. Every $k$-dimensional totally
geodesic subspace of $\Mk$ is the intersection of $\Mk$
with a $(k+1)$-dimensional linear subspace of $\RR^{d+1}$.
Isometries of $\Mk$ extend to linear transformations
of the ambient space preserving $\ink{\cdot}{\cdot}$.
Since the intensity measure
$\operatorname{Vol}_{\Mk}\otimes\mu$ is invariant under
spatial isometries, we may therefore assume without
loss of generality that $L_k=\left\{x\in\Mk:x_{k+1}=\cdots=x_d=0\right\}$.
 
\begin{proof}[Proof of Theorem~\ref{thrm:intersectiondistribution}]
Let
$\tilde L_k=\{x\in\RR^{d+1}:x_{k+1}=\cdots=x_d=0\}$,
so that $L_k=\Mk\cap\tilde L_k$.
Let $\pi_k : \RR^{d+1} \to \tilde{L}_k$ denote the orthogonal projection onto $\tilde{L}_k$. Because the coordinates $x_{k+1}, \dots, x_d$ have a positive signature $+1$ in both the Euclidean and Lorentzian inner products, $\tilde{L}_k$ is orthogonal to the kernel of $\pi_k$. Consequently, for any point $x \in \RR^{d+1}$ and any $y \in \tilde{L}_k$, we have the property
\begin{equation}\label{eq:inkprojection}
\ink{x}{y} = \ink{\pi_k(x)}{y}.
\end{equation}

Recall that the Laguerre tessellation $\mathcal{L}_\kappa(\eta_\mu)$ can be constructed by radially projecting the facets of the generalised polar region $PR(\eta_\mu) \subset \RR^{d+1}$ onto $\Mk$. The generalised polar region is defined as the set theoretic difference of $P^+(\eta_\mu) \setminus \mathrm{int} (P^{-}(\eta_\mu))$, in the negative curvature the negative set is simply the empty set. 

Both positive and negative polar intersection sets are
intersections of closed halfspaces.
When we intersect either polar intersection set with
$\tilde L_k$, we obtain the corresponding polar intersection
set in $\tilde L_k$ generated by the orthogonally projected
ambient points.
Applying \eqref{eq:inkprojection}, we obtain
\begin{align}\label{eq:positive}
P^+(\eta_\mu)\cap\tilde L_k
&=
\bigcap_{(v,t)\in\eta_\mu}
\left\{
y\in\tilde L_k:
\ink{\alpha_\kappa(v,t)}{y}\leq\kappa^{-1}
\right\}
\nonumber\\
&=
\bigcap_{(v,t)\in\eta_\mu}
\left\{
y\in\tilde L_k:
\ink{\pi_k(\alpha_\kappa(v,t))}{y}\leq\kappa^{-1}
\right\}.
\end{align}
This shows that $P^+(\eta_\mu)\cap\tilde L_k$ is precisely
the positive polar intersection set in $\tilde L_k$
generated by the projected ambient points
$\tilde x=\pi_k(\alpha_\kappa(v,t))$.
Exactly the same argument gives, for $\kappa>0$,
\begin{align}\label{eq:negative}
P^-(\eta_\mu)\cap\tilde L_k
=
\bigcap_{(v,t)\in\eta_\mu}
\left\{
y\in\tilde L_k:
\ink{\pi_k(\alpha_\kappa(v,t))}{y}\leq-\kappa^{-1}
\right\}.
\end{align}

Recall that the Laguerre cells are recovered from the
radial projections of the corresponding polar faces,
taking the closure of their union in the spherical case.
We first observe that radial projection commutes with
intersection with the subspace:
\begin{equation}\label{eq:sectionequility}
\proj_{\Mk}(F\cap\Omega_\kappa)\cap L_k
=
\proj_{\MM_\kappa^k}
\bigl(F\cap\tilde L_k\cap\Omega_\kappa\bigr).
\end{equation}
Indeed, let
$u\in\proj_{\Mk}(F\cap\Omega_\kappa)\cap L_k$.
Then $u=\proj_{\Mk}(x)$ for some
$x\in F\cap\Omega_\kappa$, so $x$ and $u$ lie on the
same positive ray through the origin.
Since $u\in L_k\subset\tilde L_k$ and $\tilde L_k$
is a linear subspace, this ray is contained in
$\tilde L_k$. Hence
$x\in F\cap\tilde L_k\cap\Omega_\kappa$.
The converse inclusion follows immediately.
 
Let $F'$ be a $k$-dimensional facet of
$P^\pm(\eta_\mu)\cap\tilde L_k$ such that
$F'\cap\Omega_\kappa\neq\emptyset$.
Condition \textup{(P1)} ensures that, near each point
of $F'\cap\Omega_\kappa$, the intersection is determined
by finitely many of its defining halfspaces.
Consequently, there exists $(v,t)\in\eta_\mu$ such that
$F'$ lies in the polar hyperplane determined by
$\pi_k(\alpha_\kappa(v,t))$.
By \eqref{eq:inkprojection}, this hyperplane is precisely
the intersection with $\tilde L_k$ of the corresponding
polar hyperplane determined by $\alpha_\kappa(v,t)$.
Since $F'$ meets $\Omega_\kappa$, the corresponding face
of $P^\pm(\eta_\mu)$ also meets $\Omega_\kappa$ and,
by general position, is almost surely a $d$-dimensional
facet $F$ of $P^\pm(\eta_\mu)$. Hence $F'=F\cap\tilde L_k$.

Conversely, let $F$ be a $d$-dimensional facet of
$P^\pm(\eta_\mu)$ such that
$F\cap\tilde L_k\cap\Omega_\kappa\neq\emptyset$.
By \eqref{eq:positive} and \eqref{eq:negative}, the intersection
$F\cap\tilde L_k$ is a face of the corresponding
polar intersection set in $\tilde L_k$.
A priori, its dimension may be smaller than $k$.
For each weighted point $(v,t)\in\eta_\mu$, we project
the vector $tv$ onto $\tilde L_k$, obtaining $\pi_k(tv)$.
The intensity measure of these vectors has a density
with respect to Lebesgue measure on $\tilde L_k$,
as the calculation below shows, even when $\mu$ has atoms.
Consequently, the vectors $\pi_k(tv)$ are almost surely
in affine general position.
Together with the finite local description
provided by \textup{(P1)}, this implies that
$\dim(F\cap\tilde L_k)=k$ almost surely.
Thus $F\cap\tilde L_k$ is a facet of the corresponding
polar intersection set in $\tilde L_k$.

In the spherical case, let
$A\subset\Mk\times(0,\infty)$ be an at most countable
set satisfying \textup{(P1)}, and set
$\tilde A:=\{\pi_k(\alpha_\kappa(v,t)):(v,t)\in A\}$.
For every $u\in L_k$, \eqref{eq:inkprojection} gives $m_A(u)=m_{\tilde A}(u)$,
where
$m_{\tilde A}(u):=\max_{ \tilde x\in\tilde A}  \ink{\tilde x}{u} $.
 
Therefore, almost surely, the facets of the positive
and negative polar intersection sets in $\tilde L_k$
that meet $\Omega_\kappa$ are precisely the intersections
$F\cap\tilde L_k$ of the corresponding facets $F$ in
$\RR^{d+1}$ for which
$F\cap\tilde L_k\cap\Omega_\kappa\neq\emptyset$.
Combining this with \cref{eq:sectionequility}, we apply
Proposition~\ref{prop:project1} in the hyperbolic case
and Theorem~\ref{prop:cellsprojection} in the spherical
case within $L_k$.
In the spherical case, the closure is taken in $L_k$,
and we use
$m_{\eta_\mu}(u)=m_{\tilde A}(u)$ for $u\in L_k$.
We conclude that $\mathcal L_\kappa(\eta_\mu)\cap L_k$
is determined by
$\tilde A=\{\pi_k(tv):(v,t)\in\eta_\mu\}$.

We now identify these projected ambient points with weighted points in $\MM_\kappa^k\times(0,\infty)$.
Each valid projected point $\tilde{x} \in \tilde{L}_k$ is mapped back to a newly weighted point $(\tilde{v}, \tilde{t}) \in \MM_\kappa^k \times (0,\infty)$. Note that for $\kappa>0$, we must assume that $\tilde{x}\neq 0$, which causes no problem in our setting because $\pi_k(\alpha_\kappa(v,t))\neq 0$ almost surely for every $(v,t)\in\eta_\mu$, since $\{ v\in\Mk : \pi_k(v)=0\}$ is a lower-dimensional subsphere and thus a measure zero set.

We want to describe the inverse mapping $(\alpha_\kappa^{(k)})^{-1}(\tilde{x}) = (\tilde{v}, \tilde{t}\,)$ such that $\alpha^{(k)}_\kappa(\tilde{v},\tilde{t}\, ) = \tilde{t} \tilde{v} = \tilde{x}$. Taking the $\kappa$-bilinear form of $\tilde{x}$ with itself yields $\ink{\tilde{x}}{\tilde{x}} = \ink{\tilde{t} \tilde{v}}{\tilde{t} \tilde{v}} = \tilde{t}^2 \ink{\tilde{v}}{\tilde{v}} = \tilde{t}^2 \kappa^{-1}$.
Solving for the new weight and spatial coordinate, we obtain the explicit inverse map
$$
(\alpha_\kappa^{(k)})^{-1}(\tilde{x}) =
\left(
\frac{\tilde{x}}{\sqrt{\kappa \ink{\tilde{x}}{\tilde{x}}}}, \,
  \sqrt{\kappa \ink{\tilde{x}}{\tilde{x}}}  
\right).
$$

We now define the transformation that maps the weighted points $(v,t)$ to their new coordinates $(\tilde{v}, \tilde{t})$. This mapping is obtained by projecting the ambient point onto $\tilde{L}_k$ and subsequently applying the inverse map $(\alpha_\kappa^{(k)})^{-1}$ derived above. The explicit transformation is
$$
T(v,t) =
\left(
\frac{\pi_k(v)}{\sqrt{\kappa \ink{\pi_k(v)}{\pi_k(v)}}}, \,
t \sqrt{\kappa \ink{\pi_k(v)}{\pi_k(v)}} 
\right).
$$
For $\kappa<0$, it is defined for all $(v,t)\in \Mk\times (0,\infty)$, while for $\kappa>0$, it is defined for all $(v,t)\in \Mk\times (0,\infty)$ such that $\pi_k(v)\neq 0$. 
 
Since the original points $(v,t)$ form a Poisson point process $\eta_\mu$ on $\MM_\kappa^d \times (0,\infty)$, their geometric transformation under $T$ yields a new point process on $\MM_\kappa^k \times (0,\infty)$. By the mapping theorem for Poisson point processes \cite[Theorem~24.16]{Klenke}, this mapped image is almost surely a Poisson point process, which we denote as $\eta_{\tilde{\mu}}$. Important note is that on $\kappa>0$ case the transformation $T$ is applied in a full measure subset of $\Mk\times(0,\infty)$.

To compute its intensity measure $\tilde{\Lambda}$, we use generalised tubular coordinates around $L_k$. A point $v \in \MM_\kappa^d$ is parametrized by its metric projection $\tilde{v} \in \MM_\kappa^k$, its orthogonal geodesic distance $r = \dk(v, L_k)$, and an angular direction $\theta \in \SPH^{d-k-1}$. By the equidistant decomposition \cite[Chapter 3.1]{Solodovnikov}, the volume element is expressed as 
$$
\Vold (\dd v) = \cos_\kappa(r)^k \sin_\kappa(r)^{d-k-1} \,\dd r \,\dd \theta \,\Volk (\dd \tilde{v}),
$$
where $\dd\theta$ denotes surface measure on
$\SPH^{d-k-1}$.
In these coordinates, the ambient projection is $\pi_k(v) = \cos_\kappa(r) \tilde{v}$. Since $\tilde{v} \in \MM_\kappa^k$, we have $\ink{\tilde{v}}{\tilde{v}} = \kappa^{-1}$. Thus, the $\kappa$-bilinear form is $\kappa \ink{\pi_k(v)}{\pi_k(v)} = \kappa \cos_\kappa(r)^2 \ink{\tilde{v}}{\tilde{v}} = \cos_\kappa(r)^2$.
Substituting this into $T(v,t)$, the spatial coordinate maps trivially to $\tilde{v}$, and the weight simplifies to $\tilde{t} = t \cos_\kappa(r)$.

By the mapping theorem, the intensity measure
$\tilde\Lambda$ evaluated on a measurable rectangle
$B\times(s_1,s_2)\subset\MM_\kappa^k\times(0,\infty)$
is given by
$$
\tilde\Lambda(B\times(s_1,s_2))
=
\int_{\MM_\kappa^d\times(0,\infty)}
\mathbf{1}_{B\times(s_1,s_2)}(T(v,t))
\,\mu(\dd t)\,\Vold(\dd v).
$$
Substituting the tubular coordinates and the simplified
transformation $T(v,t)$, we obtain
$$
\tilde\Lambda(B\times(s_1,s_2))
=
\int_B\int_{\SPH^{d-k-1}}
\int_{I_\kappa}\int_{(0,\infty)}
\mathbf{1}_{(s_1,s_2)}(u(t,r)) \,\mu(\dd t)
\cos_\kappa(r)^k\sin_\kappa(r)^{d-k-1}
 \dd r \dd\theta \Volk(\dd\tilde v),
$$
where  $u(t,r)=t \cos_\kappa(r) $ and $I_\kappa$ is the range of the distance $r$:
$I_\kappa=[0,\infty)$ for $\kappa<0$, and
$I_\kappa=[0,\pi/(2\sqrt{\kappa}))$ for $\kappa>0$.
 
Integrating over the angular coordinate $\theta$ yields the surface area of the $(d-k-1)$-dimensional unit sphere, $\omega_{d-k} = \frac{2\pi^{(d-k)/2}}{\Gamma((d-k)/2)}$.
To compute the transformed weight measure $\tilde{\mu}$, we apply the substitution $u =t \cos_\kappa(r)$, which implies $\dd u =t (-\kappa \sin_\kappa(r)) \, \dd r $.

Since $\cos_\kappa(r)=u/t$, identity
\eqref{eq:genpythagorean} gives
$$
\sin_\kappa(r)
=
\frac{1}{t}
\sqrt{\left|\frac{t^2-u^2}{\kappa}\right|}.
$$
The generalised sine is positive in the interior of
the radial range. The endpoint $r=0$ has spatial
volume zero, so we may use the strict inequalities
$u>t$ for $\kappa<0$ and $0<u<t$ for $\kappa>0$.
For $\kappa>0$, reversing the integration limits
cancels the negative sign in the change of variables.
Substituting the expressions for $\cos_\kappa(r)$,
$\sin_\kappa(r)$ and $\dd r$, and applying Tonelli's
theorem to interchange the $t$ and $u$ integrals,
the indicator restricts $u$ to $(s_1,s_2)$. We obtain
$$
\tilde\Lambda(B\times(s_1,s_2))
= 
\int_B\int_{s_1}^{s_2}
\frac{\omega_{d-k}}{|\kappa|^{(d-k)/2}}\,u^k 
\left[
\int_{J_\kappa(u)}
|t^2-u^2|^{\frac{d-k}{2}-1}
t^{1-d}\,\mu(\dd t)
\right]
\,\dd u\,\Volk(\dd\tilde v),
$$
where $J_\kappa(u)=(0,u)$ for $\kappa<0$, and
$J_\kappa(u)=(u,\infty)$ for $\kappa>0$.

The preceding calculation shows that
$\tilde\Lambda=\Volk\otimes\tilde\mu$, where
$$
\tilde\mu(\dd u)
=
\left[
\frac{\omega_{d-k}}{|\kappa|^{(d-k)/2}}\,u^k
\int_{J_\kappa(u)}
|t^2-u^2|^{(d-k)/2-1}
t^{1-d}\,\mu(\dd t)
\right]\dd u.
$$
Integrating the kernel over a compact interval in $u$
and using \eqref{eq:IC} shows that $\tilde\mu$ is locally finite.
\end{proof}

\addtocontents{toc}{\protect\setcounter{tocdepth}{0}}
\section*{Acknowledgments}
The author acknowledges financial support from the FIS 2 project “Understanding pattern formation in nature via complex analysis” CUP B53C24009510001, as well as from the Deutsche Forschungsgemeinschaft (DFG, German Research Foundation) under Germany’s Excellence Strategy (EXC 2044/2 – 390685587, Mathematics Münster: Dynamics–Geometry–Structure) and through the Priority Programme SPP 2265 “Random Geometric Systems.”
\section*{Declarations}

\subsection*{Statement on the use of generative AI}
The author finalized most of the mathematical statements and proofs by April 2026 without using generative AI models. The figures were produced using R, with code written with the assistance of ChatGPT and subsequently corrected and adjusted by the author. This work builds on \cite{AnnaPan} and was therefore published after it. As \cite{AnnaPan} underwent revisions, the present work was not finalized until November 2026. During the final revision process, the author used ChatGPT to improve the clarity and wording of the text, identify potential mathematical issues, and correct typographical errors. The author takes full responsibility for the contents of the paper.

\subsection*{Conflict of interest statement}
The author declares that they have no conflict of interest.

\subsection*{Data availability statement}
We do not analyze or generate any datasets.

\footnotesize{
\bibliographystyle{abbrv}
\bibliography{biblio}
}
 
\end{document}